\documentclass[12pt]{article}
\usepackage{amsthm,amssymb,amsmath,amsfonts}
\usepackage{graphicx}
\usepackage{hhline,array}
\usepackage{color}
\numberwithin{equation}{section}

\newtheorem{Lemma}{\bf Lemma }[section]
\newtheorem{Th}{\bf Theorem}[section]
\newtheorem{Defin}{\bf Definition}[section]
\newtheorem{Remark}{\bf Remark}[section]

\begin{document}
\sloppy

\begin{center}
\textsc{Numerical Methods of Optimal Accuracy for Weakly Singular Volterra Integral Equations}

\textbf{I.V. Boykov, A.N. Tynda}

\emph{Penza State University, 40 Krasnaya St., 440026 Penza, Russia}

E-mails:  \textcolor{blue}{\underline{boikov@pnzgu.ru}},
          \textcolor{blue}{\underline{tynda@pnzgu.ru}}.

\end{center}

\textbf{Abstract.} Weakly singular Volterra integral equations of the different types are considered. The construction of accuracy-optimal
numerical methods for one-dimensional and multidimensional equations is discussed. Since this question is closely related with the optimal
approximation problem, the orders of the Babenko and Kolmogorov \(n-\)widths of compact sets from some classes of functions have been
evaluated. In conclusion we adduce some numerical illustrations for 2-D Volterra equations.

\vspace{5mm}


\vspace{5mm}


\section{Introduction. Definitions and auxiliary statements}

Volterra integral equations have numerous applications in economy,
ecology, medicine \cite{Baker,Brunner-2004}. For a detailed study of
approximate methods for Volterra integral equations including
Abel-Volterra equations we refer to, e.g.,
\cite{Baker,Brunner-2004,Brunner-Pedas,Diogo-McKee,Tynda-18,Tynda-31,Verlan}
and the references therein.

In this paper we summarize our results concerning the construction
of numerical methods of optimal accuracy for multidimensional weakly
singular Volterra integral equations (VIEs). Some of these results
are published for the first time.

The paper is organized as follows. In Section 1, we introduce the
classes of functions being used and prove some statements concerning
the smoothness of exact solutions of VIEs. Section 2 is dedicated to
evaluation of the Babenko and Kolmogorov widths of compact sets from
introduced classes of functions. There are also constructed the
special local splines realizing the optimal estimates. In Section 3,
we describe the projective method for multidimensional VIEs based on
the approximation of the exact solutions by these splines. The
numerical example for 2-D VIE is given in Section 4.

\subsection{Classes of functions}
\begin{Defin}
  Let $Q_{r,\gamma}^{*}(\Omega,M),$ $\Omega=[0,T]^l,\;l=1,2,\ldots,$
be a class of functions $f(t_1,\ldots,t_l)$ defined on $\Omega$
and satisfying the following conditions:
\[
  \left|\frac{\partial^{|v|}f(t_1,\ldots,t_l)}
       {\partial t_1^{v_1}\cdots\partial t_l^{v_l}}\right|
  \leqslant M,\; 0\leqslant |v|\leqslant r,
\]
\[
    \left|\frac{\partial^{|v|}f(t_1,\ldots,t_l)}
       {\partial t_1^{v_1}\cdots\partial t_l^{v_l}}\right|\leqslant
    \frac{M}{(\rho(t,\Gamma_0))^{|v|-r-\zeta}},
    \; r<|v|\leqslant s,\;\; t\in\{\Omega\setminus\Gamma_0\},
\]
where \(M\) is some constant, \(0<M<\infty\), $t=(t_1,\ldots,t_l);\; v=(v_1,\ldots,v_l),\; |v|=v_1+\cdots+v_l;$ $s=r+\gamma,\:
\zeta=0,$ if $\gamma$ is an integer; $s=r+[\gamma]+1,$ $\gamma=[\gamma]+\mu,\; 0<\mu<1,\; \zeta=1-\mu,$ if $\gamma$ is
non-integer; $\Gamma_0$ is an intersection of the boundary $\Gamma$ of domain $\Omega$ with the union of coordinate planes;
\(\rho(t,\Gamma_0)=\min\limits_{i}|t_i|;\;\;\)
\end{Defin}
\begin{Remark}
    In one-dimensional case $(l=1)$ $\Gamma_0$ is defined as the point $t=0$.
\end{Remark}
\begin{Defin}
  Let $Q_{r,\gamma}^{**}(\Omega,M)$, $\Omega=[0,T]^l,\;l=1,2,\ldots,$
be a class of functions $f(t_1,\ldots,t_l)$ defined on $\Omega$
and satisfying the following conditions:
\[
  \left|\frac{\partial^{|v|}f(t_1,\ldots,t_l)}
       {\partial t_1^{v_1}\cdots\partial t_l^{v_l}}\right|\leqslant M,\; 0\leqslant |v|\leqslant
  r,
\]
\[
  \left|\frac{\partial^{|v|}f(t_1,\ldots,t_l)}
       {\partial t_1^{v_1}\cdots\partial t_l^{v_l}}\right|\leqslant \frac{M}{(\rho(t,0))^{|v|-r-\zeta}},
  \; r<|v|\leqslant s,\;\;t\ne 0,
\]
where \(M\) is some constant, $t=(t_1,\ldots,t_l);\; v=(v_1,\ldots,v_l),\; |v|=v_1+\cdots+v_l;$ $s=r+\gamma,\; \zeta=0,$ if
$\gamma$ is an integer; $s=r+[\gamma]+1,$ $\gamma=[\gamma]+\mu,\; 0<\mu<1,\; \zeta=1-\mu,$ if $\gamma$ is non-integer;
$\rho(t,0)=\min\limits_{k=1,\ldots,l}|t_k|$.
\end{Defin}
\begin{Remark}
   It is obvious that in one-dimensional case $(l=1)$ classes
   $Q_{r,\gamma}^{*}(\Omega,M)$ and $Q_{r,\gamma}^{**}(\Omega,M)$ are
   equal.
\end{Remark}
\begin{Defin}
   Let $\Omega=[0,T]^l,\;l=1,2,\ldots,\; r=1,2,\ldots,\;
0<\gamma\leqslant 1.$ Function $f(t_1,\ldots,t_l)$ belongs to the class $B_{r,\gamma}^{*}(\Omega,A)$ if the following
inequalities hold:
\[
  |f(t_1,\ldots,t_l)|\leqslant A,
\]
\[
  \left|\frac{\partial^{|v|}f(t_1,\ldots,t_l)}
       {\partial t_1^{v_1}\cdots\partial t_l^{v_l}}\right|
       \leqslant A^{|v|}|v|^{|v|}, \; 0<|v|\leqslant
       r,\;\;t\in\Omega,
\]
\[
 \left|\frac{\partial^{|v|}f(t_1,\ldots,t_l)}
      {\partial t_1^{v_1}\cdots\partial t_l^{v_l}}\right|
      \leqslant
      \frac{A^{|v|}|v|^{|v|}}{(\rho(t,\Gamma_0))^{|v|-r-1+\gamma}},
      \;\; t\in\{\Omega\setminus\Gamma_0\},\;\;r<|v|<\infty,
\]
where $A$ is a constant independent of $|v|$.
\end{Defin}
\begin{Defin}
   Let $\Omega=[0,T]^l,\;l=1,2,\ldots,\; r=1,2,\ldots,\;
0<\gamma\leqslant 1.$ Function $f(t_1,\ldots,t_l)$ belongs to the class $B_{r,\gamma}^{**}(\Omega,A)$ if the following
inequalities hold:
\[
  |f(t_1,\ldots,t_l)|\leqslant A,
\]
 \[
  \left|\frac{\partial^{|v|}f(t_1,\ldots,t_l)}
       {\partial t_1^{v_1}\cdots\partial t_l^{v_l}}\right|
       \leqslant A^{|v|}|v|^{|v|},\; 0< |v|\leqslant
       r,
\]
\[
 \left|\frac{\partial^{|v|}f(t_1,\ldots,t_l)}
      {\partial t_1^{v_1}\cdots\partial t_l^{v_l}}\right|
      \leqslant
      \frac{A^{|v|}|v|^{|v|}}{(\rho(t,0))^{|v|-r-1+\gamma}},
      \;\;t\ne 0,\;\;r<|v|<\infty,
\]
where $A$ is a constant independent of $|v|$.
\end{Defin}
\begin{Remark}
    In one-dimensional case $(l=1)$ classes $B_{r,\gamma}^{*}(\Omega,A)$
    and $B_{r,\gamma}^{**}(\Omega,A)$ are also equal.
\end{Remark}

\subsection{Smoothness of solutions}
For a simplicity of the presentation let us consider the
one-dimensional integral equation
\begin{equation}\label{SmQ1}
   x(t)-\int\limits_0^t H(t,\tau)x(\tau)d\tau=f(t),
  \;0\leqslant t\leqslant T.
\end{equation}
We introduce the following lemmas concerning the smoothness of an
exact solution of \eqref{SmQ1}:
\begin{Lemma}

\begin{enumerate}

  \item Let \(H(t,\tau)\in Q_{r,\gamma}^*([0,T],1)\) with respect to each
  variable separately, \(f(t)\in Q_{r,\gamma}^*([0,T],1)\).
  Then the unique solution \(x(t)\) of equation \eqref{SmQ1}
  belongs to the class \(Q_{r,\gamma}^*([0,T],M)\).

  \item Let \(H(t,\tau)\in B_{r,\gamma}^*([0,T],A_1)\) with respect to each
  variable separately, \(f(t)\in B_{r,\gamma}^*([0,T],A_2)\).
  Then there exist a unique solution \(x(t)\) of equation \eqref{SmQ1}
  such that \(x(t)\in B_{r,\gamma}^*([0,T],A)\).
\end{enumerate}
\end{Lemma}
\textbf{Proof.} It is well-known that the exact solution of \eqref{SmQ1} is a continuous function. Let
\(C_0=\max\limits_{t\in[0,T]}|x(t)|\), \(F_k=\max\limits_{t\in[0,T]}|f^{(k)}(t)|\),
\(H_k=\max\limits_{(t,\tau)\in[0,T]^2}\left|\frac{\partial^{k}H(t,\tau)}{\partial t^{k_1}\partial\tau^{k_2}}\right|\),
\(0\leqslant k_1,k_2\leqslant k\), \(k_1+k_2=k\), \(1\leqslant k\leqslant r\). Let also \(M_k\) be constants depending only on
the order \(k\).

In order to prove the assertions of this lemma we differentiate
formally the expression \( x(t)=f(t)+\int\limits_0^t
H(t,\tau)x(\tau)d\tau\):
\begin{equation}\label{SmQ2}
   x'(t)=f'(t)+\int\limits_0^t H'_t(t,\tau)x(\tau)d\tau+H(t,t)x(t).
\end{equation}
Since the right hand of this formula is a continuous function, there
exists the continuous derivative \(x'(t)\) of the exact solution
\(x(t)\) which can be estimated as follows \(|x'(t)|\leqslant
F_1+H_0C_0+H_1C_0T=M_1\).

Differentiating \eqref{SmQ2} one more time, we have
\[
   x''(t)=f''(t)+(H'_t+H'_{\tau})x(t)+H(t,t)x'(t)+
   \int\limits_0^t H''_{tt}(t,\tau)x(\tau)d\tau+H'_t(t,t)x(t).
\]
Taking into account that
\[
  H_k,F_k= \left\{
                  \begin{array}{ll}
                    O(1), & k\leqslant r; \\
                    O\left(\frac{1}{t^{k-r-\zeta}}\right), & k>r,
                  \end{array}
                \right.
\]
and using the previous estimate for \(|x'(t)|\) we obtain
\[
   |x''(t)|\leqslant F_2+2H_1C_0+H_0(F_1+H_0C_0+H_1C_0T)+H_2C_0T+
   H_1C_0=
\]
\[
   =F_2+H_2C_0T+3H_1C_0+H_0(F_1+H_0C_0+H_1C_0T)=M_2.
\]
Further,
\[
   |x'''(t)|\leqslant F_3+7H_2C_0+5H_1|x'|+H_0|x''|+H_3C_0T
   =M_3.
\]

Thus, differentiating \eqref{SmQ1} \(k\)-times (\(k\leqslant s\)) we
can conclude that there exists \(x^{(k)}(t)\) and the estimate holds
\[
    |x^{(k)}(t)|\leqslant
\left\{
  \begin{array}{ll}
    M_k, & k\leqslant r; \\
    \frac{M_k}{t^{k-r-\zeta}}, & r<k\leqslant s.
  \end{array}
\right.
\]

Therefore \(x(t)\in Q_{r,\gamma}^*([0,T],M)\) with \(M=\max\limits_{k=\overline{1,s}}\{M_k,1\}\) and the first statement of Lemma
is proved.

Consider now the class \(B_{r,\gamma}^*([0,T])\). It's easy to see that \(|x^{(r)}(t)|\) is bounded. The estimate for derivative
of order \(v,\; r<v<\infty,\) has the following form
\begin{equation}\label{SmQ3}
   |x^{(v)}(t)|\leqslant \frac{C_1^vv^v}{t^{v-r-1+\gamma}}+
   \frac{C_2^vv^vC_0T}{t^{v-r-1+\gamma}}+
      \frac{4^vC_3^vv^v}{t^{v-r-2+\gamma}}\leqslant
      \frac{A^vv^v}{t^{v-r-1+\gamma}},
\end{equation}
where the constants \(C_1,C_2,C_3\) depend only on \(A_1,A_2\).

Hence the exact solution \(x(t)\in B_{r,\gamma}^*([0,T],A)\).

The statements concerning the smoothness of the exact solution of
weakly singular VIEs are distributed to a case of the
multidimensional integral equations of the following form
\begin{equation}\label{Multidim-1}
     (I-K)x\equiv x(t)-\int\limits_{0}^{t_l}\cdots\int\limits_{0}^{t_1}
     h(t,\tau)g(t-\tau)x(\tau)d\tau=f(t),
\end{equation}
where \(t=(t_1,\ldots,t_l),\;\tau=(\tau_1,\ldots,\tau_l),\)
\(0\leqslant t_1,\ldots,t_l\leqslant T\);  weakly singular kernels
\(g(t-\tau)\) may have the form
\begin{equation}\label{Multidim-2}
   g(t_1,\ldots,t_l)=t_1^{r+\alpha}\cdots t_l^{r+\alpha}
\end{equation}
or
\begin{equation}\label{Multidim-3}
   g(t_1,\ldots,t_l)=(t_1^2+\cdots+t_l^2)^{r+\alpha}.
\end{equation}

Applying the analogous technique with much more complicated
computations, we have:
\begin{enumerate}
   \item If \(g(t_1,\ldots,t_l)\) has the form \eqref{Multidim-2}, \(h(t,\tau)\) has the
         continuous partial derivatives up to order \(s=r+[\gamma]+1\) and \(f(t)\in
         Q_{r,\gamma}^*(\Omega,1)\), then the exact solution \(x(t)\)
         of equation
         \eqref{Multidim-1} belongs to \(Q_{r,\gamma}^*(\Omega,M)\) with
         \(\gamma=s-r-\alpha\). If \(h(t,\tau)\) is an analytical function,
         \(f(t)\in B_{r,\gamma}^{*}(\Omega,1)\), then \(x(t)\in B_{r,\gamma}^{*}(\Omega,A)\).
   \item If \(g(t_1,\ldots,t_l)\) has the form \eqref{Multidim-3},
         \(h(t,\tau)\) has the continuous partial derivatives up to order \(s=r+[\gamma]+1\)
         and \(f(t)\in Q_{r,\gamma}^{**}(\Omega,1)\), then the exact solution
         \(x(t)\) of equation \eqref{Multidim-1} belongs to
         \(Q_{r,\gamma}^{**}(\Omega,M)\) with \(\gamma=s-r-\alpha\).
         If \(h(t,\tau)\) is an analytical function, \(f(t)\in B_{r,\gamma}^{**}(\Omega,1)\), then \(x(t)\in
         B_{r,\gamma}^{**}(\Omega,A)\).
\end{enumerate}
Note that the smoothness properties of the exact solutions of
multidimensional weakly singular Fredholm integral equations have
been investigated in the paper \cite{Vainikko2} by G. Vainikko.

%
%
\section{The optimal reconstruction of functions from
         \(Q_{r,\gamma}^{*}(\Omega,M)\), \(Q_{r,\gamma}^{**}(\Omega,M)\),
         \(B_{r,\gamma}^{*}(\Omega,A)\), and \(B_{r,\gamma}^{**}(\Omega,A)\)}

In order to construct methods of optimal accuracy for numerical solution of VIEs we need in optimal methods for the approximation
of functions from the classes \(Q_{r,\gamma}^{*}(\Omega,M)\), \(Q_{r,\gamma}^{**}(\Omega,M)\), \(B_{r,\gamma}^{*}(\Omega,A)\),
and \(B_{r,\gamma}^{**}(\Omega,A)\).

For this purpose the Babenko and Kolmogorov n-widths of compact sets
from these classes are evaluated and local splines are constructed.
The error orders of these splines coincide with the magnitudes of
the widths. The obtained assertions are diffusion of results of the
papers \cite{Boikov1,Boikov2,Boikov3,Boikov-Tynda6}.

Let us recall definitions of the Babenko and Kolmogorov
\(n-\)widths.

Let $B$ be a Banach space, $X \subset B$ be a compact set, and $\Pi:
X \to \bar X$ be a mapping of $X \subset B$ onto a
finite-dimensional space $\bar X.$
\begin{Defin}\cite{Babenko1}
    Let $L^n$ be  $n$-dimensional subspaces of the linear space $B$.
    The Kolmogorov $n$-width $d_n(X,B)$ is defined by
\begin{equation}\label{Widths-d1}
   d_n(X,B) = \inf\limits_{L^n}\sup\limits_{x \in X}\inf\limits_{u
   \in L^n} \|x-u\|,
\end{equation}
where the external infimum is calculated over all $n$-dimensional
subspaces of $L^n.$
\end{Defin}
\begin{Defin}\cite{Babenko1}
  The Babenko $n$-width $\delta_{n}(X)$ is defined by the expression
\[
  \delta_{n}(X)=\inf\limits_{\Pi:X \to R^{n}}
  \sup\limits_{x \in X}\mathrm{diam\:} \Pi^{-1} \Pi(x),
\]
where the infimum is calculated over all continuous mappings
$\Pi:X\to R^{n}$.
\end{Defin}
If the infimum in \eqref{Widths-d1} is attained for some $L^n$, this
subspace is called an extremal subspace.

The widths play the important role in the numerical analysis and
approximation theory since they have close relations to many optimal
problems such as $\varepsilon-$complexity of integration and
approximation, optimal differentiation, and optimal approximation of
solutions for the operator equations.

A detailed study of these problems in view of general theory of
optimal algorithms is given in \cite{Traub-1980}.

Throughout this paper \(A\) and \(A_k,\;k=1,2,\ldots,\) denote some positive constants that do not depend on \(N\).
\begin{Th}\label{Widths-T1}
Let $\Omega=[0,T]^l,\; l=1,2,\ldots.$ Then the estimates hold
\[       d_n(Q_{r,\gamma}^{*}(\Omega,M),C(\Omega))\asymp
       \delta_n(Q_{r,\gamma}^{*}(\Omega,M),C(\Omega))\asymp \varepsilon_n,\]
where \(\varepsilon_n \asymp n^{-s}\) if \(l=1\) and
   \begin{equation}\label{Q-estimate}
\varepsilon_n \asymp
        \begin{cases}
           n^{-(s-\gamma)/(l-1)}, & v>\frac{l}{l-1}; \\
           n^{-s/l}(\ln n)^{s/l}, & v=\frac{l}{l-1}; \\
           n^{-s/l}, & v<\frac{l}{l-1},
        \end{cases}
\end{equation} if \(l\ne 1\). Here \(v=s/(s-\gamma)\).
\end{Th}

\textbf{Proof.}  In order to estimate the infimum of the Babenko widths we divide the domain $\Omega$ into parts
$\Delta_k,\;k=0,1,\ldots,N-1$. Here $\Delta_k$ denotes the set of points from $\Omega$ satisfying the inequalities
\[
  \left(\frac{k}{N}\right)^vT\leqslant\rho(t,\Gamma_0)
  \leqslant\left(\frac{k+1}{N}\right)^vT,
  \;k=0,1,\ldots,N-1,\; t=(t_1,\ldots,t_l).
\]

Let $h_k=\left(\frac{k+1}{N}\right)^vT-
           \left(\frac{k}{N}\right)^vT,\;k=0,1,\ldots,N-1.$
Each of domains $\Delta_k$ we then cover with cubes and
parallelepipeds $ \Delta_{i_1, \ldots, i_l}^k$ which edges do not
exceed $h_k$ and parallel to the coordinate axes. For each cube
$\Delta_{i_1,\ldots,i_l}^k=[a_{i_1}^k,a_{i_1+1}^k;\ldots;
a_{i_l}^k,a_{i_l+1}^k],$ $k=0,1,\ldots,N-1$ we construct the
function of the following form
\[
  \Psi_{i_1,\ldots,i_l}^k(t)=
 \left\{%
\begin{array}{ll}
    A\frac{\Bigl((t_1-a_{i_1}^k)(a_{i_1+1}^k-t_1)
  \cdots(t_l-a_{i_l}^k)(a_{i_l+1}^k-t_l)\Bigr)^s}
  {h^{s(2l-1)}_k((k+1)/N)^{v\gamma}}, & \hbox{ for } t\in \Delta_{i_1,\ldots,i_l}^k; \\
    0, & \hbox{ for } t\notin \Delta_{i_1,\ldots,i_l}^k. \\
\end{array}%
\right.
\]
The constant $A$ is chosen from the condition \(\Psi_{i_1,\ldots,i_l}^k\in Q_{r,\gamma}^*(\Omega,M)\).

We now estimate a maximum of the function
$\Psi_{i_1,\ldots,i_l}^k(t).$ It is obvious that
\[
  \max\limits_{t\in\Omega}|\Psi_{i_1,\ldots,i_l}^k(t)|\geqslant A_1h^s
  \left(\frac{N}{k+1}\right)^{v\gamma}=
  A_2\frac{(k+\theta)^{(v-1)s}}{(k+1)^{v\gamma}}
  \frac{1}{N^{v(s-\gamma)}},\;0<\theta<1.
\]

The value $v$ is chosen so that
$\max\limits_t|\Psi_{i_1,\ldots,i_l}^k(t)|$ be independent of the
number $k$. It is obvious that the condition $v=s/(s-\gamma)$ is
sufficient.

Hence,
\begin{equation}\label{Q-31}
  |\Psi_{i_1,\ldots,i_l}^k(t)|\geqslant \frac{A_1}{N^s}.
\end{equation}

Introduce the function $\Psi_{\lambda}(t)$ by the formula
\[
  \Psi_{\lambda}(t)=\sum\limits_{k=0}^{N-1}\sum\limits_{i_1,\ldots,i_l}
  \lambda_{k,i_1,\ldots,i_l}\Psi_{i_1,\ldots,i_l}^k(t),
  \; -1\leqslant\lambda_{k,i_1,\ldots,i_l}\leqslant 1.
\]
Applying the Borsuk theorem \cite{Babenko1}, we obtain \(\delta_n(Q_{r,\gamma}^*(\Omega,M),C(\Omega))\geqslant \frac {A_1}
{N^s}\), where $n$ is the number of the cubes $\Delta_{i_1, \ldots, i_l}^k$ covering $\Omega$.

It is easy to see that for \(l=2,3,\ldots,\)
\[
  n\asymp 2^l\sum_{k=0}^{N-1}\left[\frac{ T-\left(\frac{k}{N}\right)^v T }
  {h_k}\right]^{l-1}=
  2^l\sum_{k=0}^{N-1}\left[\frac{N^v-k^v}{(k+1)^v-k^v}\right]^{l-1}\asymp
\]
\begin{equation}\label{Q-32}
  \asymp N^{v(l-1)}+2^l\sum_{k=1}^{N-1}
  \left[\frac{N^v-k^v}{(k+\theta)^{v-1}}\right]^{l-1}\asymp
  \begin{cases}
     N^{v(l-1)}, & v>\frac{l}{l-1}; \\
     N^l\ln N, & v=\frac{l}{l-1}; \\
     N^l, & v<\frac{l}{l-1},
  \end{cases}
\end{equation}

Hence,
\begin{equation}\label{Q-33}
  \delta_n(Q_{r,\gamma}^{*}(\Omega,M),C(\Omega))\geqslant
  \begin{cases}
     n^{-(s-\gamma)/(l-1)}, & v>\frac{l}{l-1}; \\
     n^{-s/l}(\ln n)^{s/l}, & v=\frac{l}{l-1}; \\
     n^{-s/l}, & v<\frac{l}{l-1},
  \end{cases}
\end{equation}

Let us construct a continuous local spline realizing estimate
\eqref{Q-33}.

At first we describe \textbf{the case of \(l=1\)} in more detail:

Let $N$ and $n$ be integers such that $n=(N-1)(s-1)+s$. We introduce
the partition of the interval $[0,T]$ with grid points
\(v_k=T(\frac{k}{N})^q\), $k=0,1,\ldots,N,$ where \(q=s/(s-\gamma)\)
if \(\gamma\) is an integer or \(q=s/(s-[\gamma]-1)\) if \(\gamma\)
is a non-integer.

Denote by $\Delta_k$ the segments $\Delta_k=[v_k,v_{k+1}],\;
k=0,1,\ldots,N-1.$ Let
\[
  \xi_j^k=\frac{v_{k+1}+v_k}{2}+\frac{v_{k+1}-v_k}{2}y_j,
\]
\[
  j=1,2,\ldots,s-2;\;\xi_0^k=v_k,\;\xi_{s-1}^k=v_{k+1};
  \;k=1,2,\ldots,N-1,
\]
\[
  \xi_j^0=\frac{v_{1}+v_0}{2}+\frac{v_{1}-v_0}{2}w_j,\;
  j=1,2,\ldots,r-2;\;\xi_0^0=v_0,\;\xi_{r-1}^0=v_{1},
\]
where $y_j$ and $w_j$ are the roots of the Legendre polynomials of
degrees $s-2$ and \(r-2\) respectively.

We denote by $P_s(f,\Delta_k)$ the operator replacing the function
\(f (t),\;t\in\Delta_k\)  by the interpolation polynomial of degree
\(s-1\) for \(k=\overline{0,N-1}\) constructed at the nodes
\(\xi_j^k\).

Let then $f_N(t)$ be a local spline defined in $[0,T]$ and
composed of polynomials \(P_s(f,\Delta_k),\;k=0,1,\ldots,N-1\).

It is easy to see that
\[
  \|f(t)-f_N(t)\|_{C[\Delta_k]}\leqslant AN^{-s}, \;
   k=\overline{1,N-1}
\]

For the segment \(\Delta_0\) we have
\[
  \|f(t)-f_N(t)\|_{C[\Delta_0]}\leqslant \frac{A_1(v_1-v_0)^r}{r!}=
  A_2\left(\frac{T}{N^q}\right)^r=\frac{A_3}{N^{qr}}=AN^{-s}.
\]

Hence, \(\|f(t)-f_N(t)\|_{C[0,T]}\leqslant AN^{-s}\). Since the
general number \(n\) of the functionals using for construction of
the spline $f_N(t)$  is estimated as \(n\asymp N\), we obtain
\[
  d_n(Q_{r,\gamma}^{*}(\Omega),C)\asymp n^{-s}
\]

\textbf{Case of \(l=2,3,\ldots\).}

Here we introduce the operator
$P_s[f,[a_1,b_1;\ldots;a_l,b_l]]=P_s^{t_1}\cdots P_s^{t_l}$ where
$P_s^{t_i}$ is the interpolating polynomial operating on variable
$t_i,\;i=1,2,\ldots,l,$ and constructed as earlier in the case of
\(l=1\). The construction of the continuous local spline starts with
domain $\Delta_{N-1}$. In this domain function $f(t_1,\ldots,t_l)$
is replaced with the interpolating polynomial $P_s[f,\Delta_{N-1}]$.
In order to construct the local spline in $\Delta_{N-2}$ we cover
this domain with the cubes and parallelepipeds
$\Delta^{N-2}_{i_1,\ldots,i_l}$. Note that their edges do not exceed
$h_{N-2}$. Here the vertices of $\Delta_{N-1}$ located on the
boundary of $\Delta_{N-2}$ are also vertices of appropriate cubes
from the set $\Delta_{i_1,\ldots,i_l}^{N-2}$. In
$\Delta_{i_1,\ldots,i_l}^{N-2}$ the function $f(t_1,\ldots,t_l)$ is
approximated by interpolating polynomials
$P_s[f,\Delta_{i_1,\ldots,i_l}^{N-2}]$.

Note that we interpolate the function \(P_s[f,\Delta_{N-1}]\)
instead of \(f(t_1,\ldots,t_l)\) at
\(t\in\Bigl\{\Delta_{i_1,\ldots,i_l}^{N-2}\bigcap
\Delta_{N-1}\Bigr\}\) (the observance of a continuity condition).

In \(\Delta_k,\;k=0,1,\ldots,N-3,\) the spline is then constructed
in a similar way.

Let \(f_s^*(t_1,\ldots,t_l)\) be the spline composed of
polynomials \(P_s[f,\Delta^k_{i_1,\ldots,i_l}]\). It is clear that
\begin{equation}\label{Q-34}
  \|f(t)-f_s^*(t)\|_C\leqslant AN^{-s}
\end{equation}

From this estimate and inequality \eqref{Q-33} it follows that the
right part of relation \eqref{Q-estimate} holds. Using the
inequality \(\delta_n\leqslant 2d_n\) \cite{Babenko1}, we finish
the proof.

\begin{Th}\label{Widths-T2}
  Let $\Omega=[0,T]$. The estimates hold
  \begin{equation}\label{t2}
     \delta_n(B_{r,\gamma}^{*}(\Omega,A),C(\Omega))\leqslant
     2^{-\sqrt{n}(r+1-\gamma)},\;
     d_n(B_{r,\gamma}^{*}(\Omega,A),C(\Omega))\asymp
     2^{-\sqrt{n}(r+1-\gamma)}.
  \end{equation}
\end{Th}
\textbf{Proof.} At first we construct a continuous local spline realizing the estimate \eqref{t2}. It will allow us to obtain an
upper bound estimate for the Kolmogorov \(n-\)width \(d_n(B_{r,\gamma}^{*}(\Omega,A),C(\Omega))\).

The interval \([0,T]\) is divided into $N+1$ segments
\(\Delta_k=[v_{k},v_{k+1}],\; k=0,\ldots,N,\) with the knots
\(v_0=0,\) \(v_k=2^{k-1-N}T,\) \(k=1,\ldots,N+1\).

Let
\begin{equation}\label{Widths-nodes}
  \begin{split}
    \xi_j^k=\frac{v_{k+1}+v_k}{2}+\frac{v_{k+1}-v_k}{2}y_j,
    \hspace*{3.3cm}& \\
     j=1,2,\ldots,m_k-2;\;\xi_0^k=v_k,\;\xi_{m_k-1}^k=v_{k+1};
    \;k=0,1,\ldots,N,
  \end{split}
\end{equation}
where \(y_j\) are the roots of the first kind Chebyshev
polynomials of degree \(m_k-2\); \(m_0=r\),
\(m_k=\Bigl[\frac{10}{9}k(r+1-\gamma)AT\Bigr]+1,\)
\(k=1,2,\ldots,N\).

Denote as earlier by $P_k(f,\Delta_k)$ the operator interpolating
the function \(f(t),\;t\in\Delta_k,\) with the polynomial of
degree \(m_k-1\) constructed at the nodes \(\xi_j^k\),
\(k=\overline{0,N}\). Denote also \(h_k=v_{k+1}-v_k,\;
k=0,1,\ldots,N\).

Let then $f_N(t)$ be a local spline defined in $[0,T]$ and
composed of polynomials \(P_k(f,\Delta_k),\;k=0,1,\ldots,N\).

It is easy to see that the approximation error at \(t\in\Delta_0\) is
\begin{equation}\label{Error0}
    \|f(t)-f_N(t)\|_{C[\Delta_0]}\leqslant A_1\ln r h_0^{r+1-\gamma}=
    A_2\left(\frac{T}{2^{N}}\right)^{r+1-\gamma}=\frac{A_3}{2^{N(r+1-\gamma)}}.
\end{equation}

Since the degrees \(m_k-1\) of interpolating polynomials increase
proportionally to the number \(k\) of a segment, it is necessary
to estimate the Lebesgue constant \(\lambda_{m_k}\) for the nodes
\eqref{Widths-nodes}. This constant is required to define the
approximation error on the segments
\(\Delta_k,\;k=\overline{1,N}\).

It is well known that the Lebesgue constant does not depend on the
length of a segment. It depends only on a distribution of the
nodes in this segment. Therefore for a simplicity of designations
let us consider the interval \([-a,a]\) and the nodes
\[
  t_j=a y_j,\;
  j=1,2,\ldots,m_k-2;\;t_0=-a,\;t_{m_k-1}=a,
\]
where \(y_j\) are the roots of the first kind Chebyshev
polynomials of degree \(m_k-2\).

Then
\[
  \lambda_{m_k}=\max\limits_{t\in[-a,a]}
  \sum\limits_{i=0}^{m_k-1}\left|\psi_{m_k-1,i}(t)\right|,
\]
where \(\psi_{m_k-1,i}(t)=\prod\limits_{j=0,j\ne
i}^{m_k-1}\frac{(t-t_j)}{(t_i-t_j)}\) are the fundamental
polynomials.

Let \(i\ne 0,\;i\ne m_k-1\). Then
\[
  |\psi_{m_k-1,i}(t)|=\left|\Phi_{m_k-1,i}(t)\varphi_{m_k-1,i}(t)\right|
  \leqslant\left|\Phi_{m_k-1,i}(t)\right|
  \max\limits_{t\in[-a,a]}|\varphi_{m_k-1,i}(t)|,
\]
where
\[
  \Phi_{m_k-1,i}(t)=\frac{(t-t_1)\cdots(t-t_{i-1})(t-t_{i+1})\cdots(t-t_{m_k-2})}
  {(t_i-t_1)\cdots(t_i-t_{i-1})(t_i-t_{i+1})\cdots(t_i-t_{m_k-2})},
\]
\[\varphi_{m_k-1,i}(t)=\frac{(t-t_0)(t-t_{m_k-1})}{(t_i-t_0)(t_i-t_{m_k-1})}.\]

It is obvious that
\[
   \max\limits_{t\in[-a,a]}|\varphi_{m_k-1,i}(t)|=
   \frac{a^2}{|(t_i+a)(t_i-a)|}
   \leqslant\frac{a^2}{a^2-t_1^2}=\frac{a^2}{a^2-a^2\cos^2\frac{1}{2(m_k-2)}}=
\]
\[
   =\frac{1}{1-\cos^2\frac{1}{2(m_k-2)}}=\frac{1}{\sin^2\frac{1}{2(m_k-2)}}
   \asymp m_k^2, \forall i.
\]

Taking into account that
\[
   \max\limits_{t\in[-a,a]}|\psi_{m_k-1,0}(t)|\asymp
   \max\limits_{t\in[-a,a]}|\psi_{m_k-1,m_k-1}(t)|\asymp m_k,
\]
we have
\[
  \lambda_{m_k}\leqslant A_4m_k^2\sum\limits_{i=1}^{m_k-2}
  \left|\Phi_{m_k-1,i}(t)\right|+|\psi_{m_k-1,0}(t)|+|\psi_{m_k-1,m_k-1}(t)|=
  O(m_k^2\ln m_k).
\]

Therefore, the approximation error
\(\|f(t)-f_N(t)\|_{C[\Delta_k]},\; k=\overline{1,N-1}\) can be
estimated as (see, e.g., \cite{Dziadyk2}):
\[
  \|f(t)-f_N(t)\|_{C[\Delta_k]}\leqslant A_5\lambda_{m_k}
  \left(\frac{v_{k+1}-v_k}{2m_k}\right)^q
  \frac{A^qq^q}{v_k^{q-r-1+\gamma}},
\]
where \(q=\left[\frac{5(r+1-\gamma)k}{9}\right]+1\) is the maximal
order of the derivatives used for the estimation of an error.

Continuing the previous inequality we have
\[
  \|f(t)-f_N(t)\|_{C[\Delta_k]}\leqslant \frac{A_6m_k^2\ln m_k
  \left(\frac{T2^{k-1-N}}{2m_k}\right)^qA^qq^q}{v_k^{q-r-1+\gamma}}=
\]
\[
  =\frac{A_7k^2\ln k T^qA^qq^q}{m_k^q 2^{(N+2-k)q}
  2^{(k-1-N)(q-r-1+\gamma)}}=
\]
\[
  =\frac{A_8k^2\ln k \;T^qA^qq^q}{\bigl(\frac{10}{9}k(r+1-\gamma)\bigr)^q
  \;A^qT^q\; 2^{N(r+1-\gamma)}\: 2^{q-(r+1-\gamma)(k-1)}}=
\]
\[
  =\frac{A_8k^2\ln k \;q^q}{\bigl(2q\bigr)^q
  \;2^{N(r+1-\gamma)}\: 2^{q-(r+1-\gamma)(k-1)}}=
\]
\[
  =\frac{A_8k^2\ln k}{2^{N(r+1-\gamma)}\; 2^{2q-(r+1-\gamma)(k-1)}}
  \leqslant \frac{A_9k^2\ln k}{2^{N(r+1-\gamma)}\;
  2^{\frac{r+1-\gamma}{9}(k+9)}}.
\]

Hence, for all sufficiently large \(N\) the estimate holds
\[
   \|f(t)-f_N(t)\|_{C[\Delta_k]}\leqslant
   \frac{A_{10}}{2^{N(r+1-\gamma)}}.
\]

Therefore for the whole segment \([0,T]\) we have
\begin{equation}\label{Error-k}
    \|f(t)-f_N(t)\|_{C[0,T]}\leqslant
     \frac{A_{11}}{2^{N(r+1-\gamma)}}.
\end{equation}

The total number \(n\) of the functionals used for the construction
of a spline can be estimated as
\[
   n=\sum\limits_{k=0}^{N}m_k=\sum\limits_{k=0}^{N}
   \frac{10}{9}k(r+1-\gamma)AT\asymp N^2.
\]

Inequality \eqref{Error-k} allows us to define an upper bound of
the Kolmogorov \(n-\)width
\[d_n(B_{r,\gamma}^{*}(\Omega,A),C(\Omega))\leqslant
\frac{A_{14}}{2^{\sqrt{n}(r+1-\gamma)}}.\]

Repeating the arguments given in \cite{Boikov3}, \cite{Boikov4} we receive 
\[
  \delta_n(B_{r,\gamma}^{*}(\Omega,A),C(\Omega))\geqslant
  \frac{A_{13}}{2^{\sqrt{n}(r+1-\gamma)}}.
\]

Taking into account the inequality \(\delta_n\leqslant 2d_n\), we
accomplish the proof.

Note that instead of \eqref{Widths-nodes} we can also use the
another system of the nodes:
\begin{equation}\label{Widths-nodes-2}
    \xi_j^k=\frac{v_{k+1}+v_k}{2}+\frac{v_{k+1}-v_k}{2}y_j,
    \; j=0,2,\ldots,m_k-1;\; \;k=0,1,\ldots,N,
\end{equation}
where \(y_j\) are the roots of the first kind Chebyshev polynomials
of degree \(m_k-1\); \(m_0=r\),
\(m_k=\Bigl[k(r+1-\gamma)AT\Bigr]+1,\) \(k=1,2,\ldots,N\).

This allows us to eliminate the additional multiplier \(m_k^2\) in
the estimate of the Lebesque constant \(\lambda_{m_k}\). However,
the closed system of nodes \eqref{Widths-nodes} is more suitable in
practice for the numerical solution of VIEs by the projective method
described in Section 3.

\begin{Th}\label{Widths-T4}
   Let $\Omega=[0,T]^l,\; l=2,3,\ldots$. Then the estimates hold
   \begin{equation}\label{t4}
      \delta_n(Q_{r,\gamma}^{**}(\Omega,M),C(\Omega))\asymp
       d_n(Q_{r,\gamma}^{**}(\Omega,M),C(\Omega))\asymp n^{-s/l}.
   \end{equation}
\end{Th}

\textbf{Proof.} In order to estimate the infimum of \(\delta_n(Q_{r,\gamma}^{**}(\Omega,M),C(\Omega))\) we cover domain $\Omega$
with cubes as follows. The cube \(\Delta^1=\Delta_{1,\ldots,1}^1\) is an intersection of domains
\[
  \Bigl(0\leqslant t_1\leqslant
  \left(\frac{1}{N}\right)^vT\Bigr)\Bigg.\cap \cdots \Bigg.\cap
  \Bigl(0\leqslant t_l\leqslant \left(\frac{1}{N}\right)^vT\Bigr),
\]
\(v=s/(s-\gamma)\) if \(\gamma\) is an integer,
\(v=s/(s-[\gamma]-1)\) if  \(\gamma\) is a non-integer.

The domain $\Delta^2$ is then defined as
\(\Delta^2=\Delta_2'\setminus \Delta_1''\), where
\[
  \Delta_k'=\left\{(t_1,\ldots,t_l):\;0\leqslant t_1,\ldots,t_l\leqslant
  \left(\frac{k}{N}\right)^vT \right\},
\]
\[
  \Delta_k''=\left\{(t_1,\ldots,t_l):\;0\leqslant t_1,\ldots,t_l<
  \left(\frac{k}{N}\right)^vT \right\}.
\]

This domain is covered with cubes and parallelepipeds
\(\Delta_{i_1,\ldots,i_l}^2\) which edges are parallel to the axes
of coordinates and do not exceed
\(h_1=\left(\frac{2}{N}\right)^vT-\left(\frac{1}{N}\right)^vT\). The
further construction is carried out by analogy.

Each domain $\Delta^k=\Delta_k'\setminus \Delta_{k-1}'',\;
k=3,\ldots,N-1,$ is covered with cubes and parallelepipeds
$\Delta^k_{i_1,\ldots,i_l}$ with edges not exceeding
$h_{k-1}=\left(\frac{k}{N}\right)^vT-\left(\frac{k-1}{N}\right)^vT$.

In $\Delta^k_{i_1,\ldots,i_l}$ the function
$\Psi_{i_1,\ldots,i_l}^k$ is defined and then in domain $\Omega$ the
function $\Psi_{\lambda}(t_1,\ldots,t_l)$ is introduced (by analogy
to the proof of theorem 2.1). Then we show that
\(|\Psi_{\lambda}(t_1,\ldots,t_l)|\geqslant \frac{A_1}{N^s}\).

Let us define the number $n$ of parallelepipeds
$\Delta^k_{i_1,\ldots,i_l}$. It is easy to see that
\[
  n\asymp\sum_{k=1}^{N-1}\left[\frac{\left(\frac{k+1}{N}\right)^v}
  {\left(\frac{k+1}{N}\right)^v-\left(\frac{k}{N}\right)^v}\right]^{l-1}
  \asymp\sum_{k=1}^{N-1}\left[\frac{(k+1)^v}{(k+\theta)^{v-1}}\right]^{l-1}
  \asymp\sum_{k=1}^{N-1}k^{l-1}\asymp N^l.
\]

Hence, $ \delta_n(Q_{r,\gamma}^{**}(\Omega,M),C(\Omega))\geqslant An^{-s/l}$. The construction of the local spline
$f_N^{**}(t_1,\ldots,t_l)$ and further argumentation are carried out by analogy to the proof of Theorem \ref{Widths-T1}. The
theorem is proved.

\begin{Th}\label{Widths-T5}
   Let $\Omega=[0,T]^l,\; l=2,3,\ldots,
   \;0<\gamma\leqslant 1$. Then the estimates hold
  \begin{equation}\label{t5}
     \delta_n(B_{r,\gamma}^{*}(\Omega,A),C(\Omega))\asymp
     d_n(B_{r,\gamma}^{*}(\Omega,A),C(\Omega))\asymp
     \frac{1}{n^{(r+1-\gamma)/(l-1)}}.
  \end{equation}
\end{Th}

\textbf{Proof.} Let $\Delta_0$ be a set of points $t\in\Omega$
such that $0\leqslant\rho(t,\Gamma_0)\leqslant 2^{-N}T,$ and
$\Delta_k,$ $k=1,2,\ldots,N,$ be a set of points $t\in\Omega$ such
that
\[
  \frac{2^{k-1}}{2^N}T\leqslant\rho(t,\Gamma_0)\leqslant \frac{2^k}{2^N}T.
\]

Let us cover each domain $\Delta_k,\;k=0,1,\ldots,N,$ with cubes
$\Delta^k_{i_1,\ldots,i_l}$. The edges of these cubes are parallel
to the edges of $\Omega$. These edges are not less than \(h_k\) and
not more than \(2h_k\), where $h_k=\frac{2^{k-1}}{2^N}T,$
$k=0,1,\ldots,N-1$.

Now we estimate a general number of elements
$\Delta^k_{i_1,\ldots,i_l}$ covering domain $\Omega$. It is
obvious that
\[
  n\asymp \sum_{k=1}^N\left[\frac{1-\frac{2^{k-1}}{2^N}}
  {\frac{2^k}{2^N}-\frac{2^{k-1}}{2^N}}\right]^{l-1}
  \asymp \sum_{k=1}^N\left(2^{N-k+1}-1\right)^{l-1}
  \asymp \sum_{k=1}^N\frac{2^{(N+1)(l-1)}}{2^{k(l-1)}}=
\]
\[
  =\frac{1}{2^{l-1}-1}\Bigl(2^{(N+1)(l-1)}-2^{l-1}\Bigr).
\]

Thus, $n\asymp 2^{N(l-1)}$.

Repeating the arguments of the paper \cite{Boikov3} we obtain the estimate
\(\delta_n(B_{r,\gamma}^{*}(\Omega,A),C(\Omega))\geqslant 2^{-N(r+1-\gamma)}\) and conclude
\[
  \delta_n(B_{r,\gamma}^{*}(\Omega,A),C(\Omega))\geqslant
  \frac{1}{n^{(r+1-\gamma)/(l-1)}}.
\]

The construction of a continuous local spline realizing estimate
\eqref{t5} is similar to construction given in the Theorem
\ref{Widths-T1}. Here the parameter \(s\) is equal to
\(s=\Bigl[\frac{10}{9}N(r+1-\gamma)AT\Bigr]+1\).

Taking into account well known inequality $\delta_n\leqslant 2d_n$
connecting the Babenko and Kolmogorov \(n-\)widths, we finish the
proof.

\begin{Th}\label{Widths-T6}
  Let $\Omega=[0,T]^l,\; l=2,3,\ldots$. Then the estimates hold
  \begin{equation}\label{t6}
     \delta_n(B_{r,\gamma}^{**}(\Omega,A),C(\Omega))\asymp
     d_n(B_{r,\gamma}^{**}(\Omega,A),C(\Omega))\asymp
     \frac{1}{2^{\sqrt[l+1]{n}\:(r+1-\gamma)}}.
  \end{equation}
\end{Th}
The proof is carried out by analogy to the proof of Theorem
\ref{Widths-T5}.

\section{Approximate solution of multidimensional VIEs}

In this section we consider the multidimensional VIEs of the form
\begin{equation}\label{M1}
     (I-K)x\equiv x(t)-\int\limits_{0}^{t_l}\cdots\int\limits_{0}^{t_1}
     h(t,\tau)g(t-\tau)x(\tau)d\tau=f(t),
\end{equation}
where \(t=(t_1,\ldots,t_l),\;\tau=(\tau_1,\ldots,\tau_l),\)
\(0\leqslant t_1,\ldots,t_l\leqslant T\); weakly singular kernel
\(g(t-\tau)\) has the form
\begin{equation}\label{M2}
   g(t_1,\ldots,t_l)=t_1^{r+\alpha}\cdots t_l^{r+\alpha}
\end{equation}
\subsection{Numerical scheme}
We look for an approximate solution of \eqref{M1} as the spline
\(x_N^*(t_1,\ldots,t_l)\) with the unknown values
\(x_N^*(\xi_{i_1}^k,\ldots,\xi_{i_l}^k),\;
(\xi_{i_1}^k,\ldots,\xi_{i_l}^k)\in
\Delta_{i_1,\ldots,i_l}^k,\;k=0,1,\ldots,N-1,\) at the knots of
the grid.

The construction of the spline \(x_N^*(t_1,\ldots,t_l)\) is
described in Section 2  and depends on the considered class of
function.

The values \(x_N^*(\xi_{i_1}^k,\ldots,\xi_{i_l}^k)\) in each cube
\(\Delta_{i_1,\ldots,i_l}^k,\;k=0,1,\ldots,N-1,\) are determined
step-by-step by the spline-collocation technique from the systems
of linear equations
\begin{equation}\label{MultiDim-4}
  \begin{split}
   (I-K)P_N[x(t),\Delta_{i_1,\ldots,i_l}^k]
   \equiv P_N[x(t),\Delta_{i_1,\ldots,i_l}^k]-\\
   -P_N\Biggl[\int\cdots\int\limits_{\hspace{-1cm}\Delta_{i_1,\ldots,i_l}^k}
   P_N^{\tau}[h(t,\tau)]g(t-\tau)
   P_N[x(\tau),\Delta_{i_1,\ldots,i_l}^k]
   d\tau,\;\Delta_{i_1,\ldots,i_l}^k \Biggr]=\\
   =P_N[f_{i_1,\ldots,i_l}^k(t),\Delta_{i_1,\ldots,i_l}^k].
  \end{split}
\end{equation}

Here $P_N$ is an operator of projection on the set of the local
splines of the form \(x_N^*(t_1,\ldots,t_l)\);
\(f_{i_1,\ldots,i_l}^k(t_1,\ldots,t_l)\) is a new right part of
equation \eqref{M1} including the integrals over domains
\(\Delta_{i_1,\ldots,i_l}^j,\;j=0,1,\ldots,k,\) processed at the
previous steps (in these domains the spline values are already
known).

All the integrals in \eqref{MultiDim-4} are calculated using the
Gauss-type cubature formulas.

There are several ways for choosing the numeration of the
subdomains \(\Delta_{i_1,\ldots,i_l}^k,\;k=0,1,\ldots,N-1\). One
of such ways allowing the parallelization of the computing process
we indicate in \cite{Boikov-Tynda3}.

\subsection{Convergence substantiation}
Let us rewrite equation \eqref{M1} and  projective method
\eqref{MultiDim-4} in the operator form:
\begin{equation}\label{MultiDim-2.1}
   x-Kx=f,\; K:X\to X,\; X\subset C(\Omega),\;
   \Omega=[0,T]^l,\; l=2,3,\ldots,
\end{equation}
\begin{equation}\label{MultiDim-2.2}
   x_N-P_NKx_N=P_Nf,\;  P_N:X\to X_N,\; X_N\subset C(\Omega),
\end{equation}
where \(X\) is one of the sets \(Q_{r,\gamma}^{*}(\Omega,M)\) or \(B_{r,\gamma}^{*}(\Omega,A)\); \(X_N\) are the sets of
corresponding local splines.

Since the homogenous Volterra integral equation \(x-Kx=0\) has
only the trivial solution, the operator \(I-K\) is injective.
Hence, the operator \(I-K\) has the bounded inverse operator
\((I-K)^{-1}:X\to X.\) For all sufficiently large \(N\) we have
the estimates
\[
  \|(I-P_NK)^{-1}\|_{C(\Omega)}=\|\Bigl((I-K)+(K-P_NK)\Bigr)^{-1}\|_{C(\Omega)} \leqslant
\]
\[
   \leqslant\frac{\|(I-K)^{-1}\|_{C(\Omega)}}{1-\|(I-K)^{-1}\|_{C(\Omega)}
   \|K-P_NK\|_{C(\Omega)}}
   \leqslant 2\|(I-K)^{-1}\|_{C(\Omega)}=A\;(const)
\]
if
\[
  \|K-P_NK\|_{C(\Omega)}\leqslant
  \frac{1}{2 \|(I-K)^{-1}\|_{C(\Omega)}}.
\]

Let us show that the last estimate holds for all sufficiently large \(N\). Since \(y(t)\equiv\bigl(Kx\bigr)(t)\in X\) and \(X\)
is a dense set in \(C(\Omega)\) (this is valid for \(Q_{r,\gamma}^{*}(\Omega,M)\) and \(B_{r,\gamma}^{*}(\Omega,A)\) ), we have
\[
  \|K-P_NK\|_{C(\Omega)}=
  \sup\limits_{x\in X,\|x\|\leqslant1}
  \max\limits_{t\in\Omega}|x(t)-P_Nx(t)|\leqslant\varepsilon_N,
\]
where \(\varepsilon_N\to 0\) as \(N\to\infty\). Therefore,
\(\|K-P_NK\|_{C(\Omega)}\leqslant\frac{1}{2\|(I-K)^{-1}\|}\)
starting with sufficiently large \(N\).

Thus, the operators \((I-P_NK)^{-1}\) are exist and uniformly
bounded and equation \eqref{MultiDim-2.2} has a unique solution
for all sufficiently large \(N\). Taking into account that
\(P_Nx\to x\) as \(N\to\infty\) for all \(x\in X\), we apply the
projection operator \(P_N\) both to the left and the right parts
of equation \eqref{MultiDim-2.1}:
\[
  x-P_NKx=P_Nf+x-P_Nx.
\]

Subtracting this equation from \eqref{MultiDim-2.2}, we obtain
\[
  (I-P_NK)(x_N-x)=P_Nx-x,
\]
\[(x_N-x)=(I-P_NK)^{-1}(P_Nx-x).\]

This implies
\begin{equation}\label{MultiDim-2.3}
    \|x_N-x\|_C\leqslant A\|P_Nx-x\|_C\leqslant \varepsilon_N(X).
\end{equation}

Thus, the accuracy of the approximate solution obtained via
projective method \eqref{MultiDim-2.2} is determined by the accuracy
\(\varepsilon_N(X)\) of the approximation of functions from \(X\) by
the local splines.

On the other hand, it follows from the Theorems given in Section 2 that for the functions from \(X\) the order of estimate
\eqref{MultiDim-2.3} cannot be improved (see Theorems 2.1 and 2.2). Hence, we conclude that algorithm \eqref{MultiDim-4} is of
optimal accuracy order on the classes \(Q_{r,\gamma}^{*}(\Omega,M)\) and \(B_{r,\gamma}^{*}(\Omega,A)\).

It was proved in paper \cite{Tynda-18} that such numerical methods
for VIEs are also optimal with respect to complexity order.

The numerical solution of multidimensional VIEs with the optimal
accuracy requires a huge number of arithmetical operations. In
\cite{Boikov-Tynda3}, employing the 2-D VIE case as an example, we
investigate the problem of accelerating the computing process by
using multiprocessor computers.
\begin{Remark}
   The suggested algorithm can be also applied to
numerical solution of multidimensional VIEs of the form
\begin{equation}\label{Mullti-2}
  \begin{split}
     x(t)-\int\limits_{0}^{t_l}\cdots\int\limits_{0}^{t_1}
     h(t,\tau)\Bigl[(t_1-\tau_1)^2+\cdots+(t_l-\tau_l)^2\Bigr]^{r+\alpha}
     x(\tau) d\tau_1\ldots d\tau_l=f(t),\\
     t=(t_1,\ldots,t_l),\;\;
     \tau=(\tau_1,\ldots,\tau_l),\;\;t\in[0,T]^l,
  \end{split}
\end{equation}
with coefficients from the classes \(Q_{r,\gamma}^{**}(\Omega,M)\) or \(B_{r,\gamma}^{**}(\Omega,A)\).
\end{Remark}

\section{Numerical illustration}

As a numerical example we consider the integral equation with coefficients from the class \(B^*_{2,0.5}(\Omega,1)\) (at the same
time these coefficients are in \(Q^*_{r,\gamma}(\Omega,1)\)  with \(r=2\) and a certain value of \(\gamma\)).
\begin{equation}\label{2D-Example-1}
   x(t_1,t_2)-\int\limits_{0}^{t_2}\int\limits_{0}^{t_1}
   (t_1-\tau_1)^{\frac{5}{2}}(t_2-\tau_2)^{\frac{5}{2}}x(\tau_1,\tau_2)d\tau_1d\tau_2
   =(t_1t_2)^{\frac{5}{2}}+\frac{25\pi^2t_1^6t_2^6}{1048576},
\end{equation}
where \((t_1,t_2)\in\Omega=[0,1]^2.\) The exact solution of
\eqref{2D-Example-1} is \(x(t_1,t_2)=(t_1t_2)^{2.5}\).

Since \(B_{r,\gamma}^*(\Omega,A)\subset Q_{r,\gamma}^*(\Omega,M)\), two different algorithms of approximation have been applied
to \eqref{2D-Example-1} for each of these classes (see their description in Section 2). The tables of results are given below.

\begin{center}
\begin{table}[ht]
{\centering
  \begin{tabular}{|c|c|c|c|c|c|c|c|}\hline
$\mathbf{N}$     &  1 &   2& 3& 5 & 10 & 15& 20\\ \hline
     $\mathbf{\varepsilon_1}$ & $7.17e-3$ & $1.12e-5$
      & $2.17e-7$ & $4.89e-8$
      & $6.13e-9$ & $6.32e-11$ & $5.89e-13$ \\ \hline
     $\mathbf{\varepsilon_2}$ & $0.011$ & $6.15e-4$ &
     $9.05e-6$ & $6.73e-7$ &
     $2.39e-8$ & $6.84e-10$ & $ 2.87e-11$
     \\  \hline
  \end{tabular}
  \caption{The error on the class \(Q_{2,2.5}^*(\Omega,1)\) }}
\end{table}
\end{center}

\begin{center}
\begin{table}[ht]
{\centering
  \begin{tabular}{|c|c|c|c|c|c|c|c|}\hline
$\mathbf{N}$     &  1 &   2& 3& 5 & 10 & 15& 20\\ \hline
     $\mathbf{\varepsilon_1}$ & $6.25e-4$ &$1.04e-7$
      & $3.10e-8$ & $5.99e-9$
      & $5.23e-10$ & $7.31e-13$ & $6.19e-15$ \\ \hline
     $\mathbf{\varepsilon_2}$ & $0.002$ & $7.45e-6$ &
     $9.73e-7$ & $7.62e-8$ &
     $1.48e-9$ & $6.35e-12$ & $3.17e-13$
     \\  \hline
  \end{tabular}
  \caption{The error on the class \(B^*_{2,0.5}(\Omega,1)\)}}
\end{table}
\end{center}

Here \(\mathbf{N}\) is the number of subdomains of the main
partition for \(\Omega\);

\(\varepsilon_1=\max\limits_{i,j}|x(t_i,t_j)-x_N(t_i,t_j)|\) is
the error at the nodes of the grid;

\(\varepsilon_2=\|x(t_1,t_2)-x_N(t_1,t_2)\|_{C(\Omega)}\) is the
error in \(\Omega\).



\begin{thebibliography}{00}
   \bibitem{Baker}
   C.T.H. Baker, A perspective on the numerical treatment
                of Volterra equations,
     {\em J. Comp. Appl. Math.} Vol.125(2000), 217-249.
  \bibitem{Brunner-2004}
  H. Brunner,
  {\em Collocation methods for Volterra integral and related
  functional differential equations}, Cambridge University Press, Cambridge, 2004.
\bibitem{Brunner-Pedas}
  H. Brunner, A. Pedas, G. Vainikko,
   The piecewise polynomial collocation method for
   nonlinear weakly singular Volterra equation,
   {\em Math. Comp.}, Vol.68, N227 (1999), 1079-1095.
  \bibitem{Diogo-McKee}
  T. Diogo, S. McKee, T. Tang,
        Collocation methods for second-kind Volterra integral equations
        with weakly singular kernels. Proc. Roy. Soc. Edin. 124A,
        1994, 199-210.
  \bibitem{Tynda-18}
   A.N. Tynda,  Numerical  algorithms of optimal complexity for weakly singular
     Volterra integral equations, {\em Comp. Meth. Appl. Math.},
     Vol.6(2006) No. 4, p.436-442.
  \bibitem{Tynda-31}
   A.N. Tynda, Spline-collocation technique for 2D weakly
               singular Volterra integral equations.
               \emph{Bulletin of Middle-Volga Math. Society}, Vol.10, 2008,
                No.2, 68-78.
  \bibitem{Verlan}
    A.F. Verlan, V.S. Sizikov \emph{Integral equations: methods, algorithms, programms.}  Kiev, Naukova Dumka,
    1986. [In Russian]
   \bibitem{Vainikko2}
    G.M. Vainikko, On the smoothness of solution of
    multidimensional weakly singular integral equations.
    {\em Mat. Sbornik.} 1989, 180, N12, 1709-1723 [In Russian]
     \bibitem{Boikov1}
     I.V. Boykov, {\em The optimal methods of approximation of the functions and computing the integrals},
     Penza, Penza State University Publishing House, 2007, 236p.[In Russian]
     \bibitem{Boikov2}
     I.V. Boikov, Approximation of Some Classes of Functions by Local Splines,
     {\em Computational Mathematics and Mathematical Physics,} Vol.38, No. 1, 1998
     21-30.
     \bibitem{Boikov3}
     I.V. Boikov,  The optimal algorithms of recovery of the functions and
     computing of the integrals on a class of infinitely
     differentiable functions, {\em Izvestia Vuzov. Matematika}, 1998, 9,14-20.
  \bibitem{Boikov-Tynda6}
  I.V. Boikov and A.N. Tynda,
  Accuracy-optimal approximate methods for solving Volterra integral
  equations, {\em Differential Equations}, Vol.38, N.9, 2002, 1305-1313.
  \bibitem{Babenko1}
   K.I. Babenko, \emph{Theoretical Foundations and Construction of Numerical Algorithms
   for Problems in Mathematical Physics}, Nauka, Moscow, 1979. [In Russian]
     \bibitem{Traub-1980}
      J.F. Traub and H. Wozniakowski, \emph{A General Theory of Optimal
     Algorithms}, Academic Press, New York, 1980.
  \bibitem{Dziadyk2}
  V.K. Dziadyk, {\em Introduction in Theory of Uniform Approximation of the Functions
   by Polynomials}.  Moscow, Nauka, 1977, 512p. [in Russian]
   \bibitem{Boikov4}
     I.V. Boykov, Optimal approximation and Kolmogorov widths estimates
for certain singular classes related to equations of mathematical physics , {\em arXiv}, math. DG/0303109.


  \bibitem{Boikov-Tynda3}
    I.V. Boikov and A.N. Tynda,
  Methods of optimal accuracy
  for approximate solution of second-kind weakly singular Volterra
  integral equations for multiprocessor computers, in Proc. of the ICCM-2002, Novosibirsk, 381-388.






\end{thebibliography}
\end{document}